\documentclass[a4paper]{eccomas_paper-2024}
\usepackage{graphicx}

\usepackage{hyperref}

\usepackage{graphicx}
\usepackage{amsmath}
\usepackage{amsfonts}
\usepackage{amssymb}
\usepackage{booktabs}
\usepackage{lipsum}  
\usepackage{xcolor}
\usepackage{amsmath,amssymb}
\usepackage{graphicx}

\DeclareMathOperator{\E}{\mathbb{E}}
\usepackage{bm}
\usepackage{subcaption}

\usepackage{setspace}
\usepackage{tikz}

\newcommand{\norm}[1]{\left\lVert#1\right\rVert}

\usepackage{etoolbox}
\patchcmd{\thebibliography}{\section*{\refname}}{}{}{}

\title{Modelling Sand Ripples in Mine Countermeasure Simulations by means of  Stochastic Optimal Control }

\author{P. Blondeel$^{1}$, F. Van Utterbeeck$^{1}$ and B. Lauwens$^{1}$}

\heading{P. Blondeel, F. Van Utterbeeck and B. Lauwens}

\address{$^{1}$ Royal Military Academy, Department of mathematics, Avenue de la Renaissance 30\\
1000 Brussels Belgium, \{philippe.blondeel, filip.vanutterbeeck, ben.lauwens\}@mil.be
}

\keywords{Stochastic Optimal Control, Mine Countermeasures, Sand Ripples}

\abstract{Modelling and simulating mine countermeasures (MCM) search missions performed by autonomous vehicles equipped with a sensor capable of detecting mines at sea is a challenging endeavour. In this work, we present a novel way to model and account for sand ripples present on the bottom of the ocean while calculating trajectories for the autonomous vehicles by means of a stochastic optimal control framework. It is known from the scientific literature that these ripples impact the sea mine detection capabilities of the autonomous vehicles. }

\begin{document}
\thispagestyle{empty}

\section{INTRODUCTION}
Modelling and simulating mine countermeasures (MCM) search missions performed by autonomous vehicles is a challenging endeavour. The goal of these simulations typically consists of computing trajectories for the autonomous vehicles in a designated zone such that the residual MCM risk in said zone is below a certain threshold while simultaneously ensuring that the mission time stays below a certain value.  This type of problem is referred to as a Coverage Path Planning (CPP) problem, see \cite{AI2021110098,HowieChoset,Abreu,6380733,GALCERAN20131258,https://doi.org/10.1049/rsn2.12256,1511771,9121686}.  From the previously cited works, two main approaches emerge. The first one consists of decomposing  the to-be surveyed zone into a grid made from square or hexagonal cells. This grid is then used to compute a trajectory for the autonomous vehicle such that each cell is visited at least once. The second approach consists of computing a trajectory according to a pattern such as, a boustrophedon pattern, a zigzag search pattern or a weaving pattern. However, a novel third approach is presented in the work of  \cite{9140316}.  There, the CPP problem is formulated as a stochastic optimal control problem. A major advantage of this approach consists of its extensibility. Meaning that from an implementation standpoint, it is straightforward for the end-user to add additional autonomous vehicles in the to-be surveyed zone. The trajectories for the multiple vehicles are then computed simultaneously resulting in a lower mission time.

From the work of \cite{8867065}, it is known that the characteristics of the ocean floor impact the detection of mines. The presence of sand ripples, which can be viewed as sand dunes on the ocean floor,   impact the detection capabilities of the autonomous vehicles when arising in the to-be surveyed zone. In order to be able to detect  mines when sand ripples are present, the autonomous vehicles need to traverse the ripples in a  perpendicular way.

In this paper, we start from the work of \cite{9140316}, and implement the MCM search mission formulation in a stochastic optimal control framework, see \cite{pulsipher2022unifying}, such that the  mission time is minimized while ensuring that the residual MCM risk in the to-be surveyed zone is below a certain threshold. The main contribution of this work consists of the novel formulation and implementation used to account for sand ripples when computing trajectories by means of a stochastic optimal control framework.

The paper is structured as follows. First we present the methodology, where we introduce the sensor model, the formulation of the stochastic optimal control problem, and our modelization of the sand ripples in the stochastic optimal control framework. Second, we present results where we show how the trajectories are adapted when sand ripples are present in the to-be surveyed zone. We present results for up to two autonomous vehicles acting in the same zone.

\section{METHODOLOGY}

In this section, we first present the equations describing the sensor model for a Forward Looking Sensor (FLS), i.e., a sensor which is only capable of detecting the presence of sea mines when they are located in front of the autonomous vehicle. Hereafter, we give the equations pertaining to the stochastic optimal control problem as set forward in \cite{9140316}. Last, we present our approach on how to account for sand ripples when using the stochastic optimal control framework.

\subsection{Sensor Model}

We briefly introduce the equations pertaining to the modelisation of a FLS. A more detailed description can be found in \cite{9140316}.

The sensor model is given by
\begin{equation}
\gamma\left(\bm{x}\left(t\right),\bm{\omega}\right) := \lambda\,p\left(\bm{x}\left(t\right),\bm{\omega}\right)\,F_\alpha\left(\bm{x}\left(t\right),\bm{\omega}\right)\,F_\varepsilon\left(\bm{x}\left(t\right),\bm{\omega}\right),
\label{eq:sensor_model}
\end{equation}
where $\lambda$ stands for the Poisson scan rate in ${s}^{-1}$.

The sensor model of Eq.\,\eqref{eq:sensor_model} consists of three parts. The first part consists of $p\left(\bm{x}\left(t\right),\bm{\omega}\right)$ and is given by
\begin{equation}
p\left(\bm{x}\left(t\right),\bm{\omega}\right) := \Phi\left(\frac{\text{FOM} - 20 \log_{10}\left(\norm{\bm{\omega} - \bm{x}\left(t\right)} + a\norm{\bm{\omega}-\bm{x}\left(t\right)}\right)}{\sigma}\right),
\end{equation}
where  $\Phi\left(\cdot\right)$ stands for the cumulative density function (CDF) of the normal distribution, $\text{FOM}$ is a parameter related to the sonar characteristics, $a$ is the attenuation coefficient in $dB/km$, and $\bm{x}\left(t\right)$ is given in Eq.\,\eqref{eq:pos}.
We give $\norm{\bm{\omega} - \bm{x}\left(t\right)}$ as
\begin{equation}
\norm{\bm{\omega} - \bm{x}\left(t\right)} := \sqrt{\left(\omega_x - x\left(t\right)\right)^2} + \sqrt{\left(\omega_y - y\left(t\right)\right)^2},
\end{equation}
where $\bm{\omega}$ is defined as  $\bm{\omega}:= \left(\omega_x,\omega_y\right)$, and stands for the position of the possible target, i.e., a sea mine.
This second part, $F_\alpha\left(\bm{x}\left(t\right),\bm{\omega}\right)$, models the detection of the sensor in front of  the autonomous vehicle according to its Field Of View (FOV), and is given as
\begin{equation}
F_\alpha\left(\bm{x}\left(t\right),\bm{\omega}\right):=\frac{1}{1+e^{p_\alpha\left(-\frac{\alpha_{\text{FOV}}}{2} - \alpha^b\left(\bm{x}\left(t\right),\bm{\omega}\right)\right)}} + \frac{1}{e^{p_\alpha\left(\alpha^b\left(\bm{x}\left(t\right),\bm{\omega}\right)-\frac{\alpha_{\text{FOV}}}{2}\right)}} - 1,
\label{eq:F_a}
\end{equation}
where $p_\alpha$ is a parameter used to adjust the slope of the sigmoidal curves, and
\begin{equation}
  \begin{gathered}
  \alpha^b\left(\bm{x}\left(t\right),\bm{\omega}\right) := \arctan2\left(dx^b\left(\bm{x}\left(t\right),\bm{\omega}\right),dy^b\left(\bm{x}\left(t\right),\bm{\omega}\right)\right) \\
dx^b\left(\bm{x}\left(t\right),\bm{\omega}\right):=\left(\omega_x -x(t)\right)\,\cos(\psi(t))+\left(\omega_y-y(t)\right)\,\sin(\psi(t)) \\
dy^b\left(\bm{x}\left(t\right),\bm{\omega}\right):=-\left(\omega_x -x(t)\right)\,\sin(\psi(t))+\left(\omega_y-y(t)\right)\,\cos(\psi(t)).
\end{gathered}
\label{eq:F_a_expanded}
\end{equation}

%
%
 In Eq.\,\eqref{eq:F_a},  $\alpha_\text{FOV}$ is the Field Of View angle of the FLS in degrees. In Eq.\,\eqref{eq:F_a_expanded}, $\psi(t)$ represents the angle the autonomous vehicle has with respect to the horizontal axis in degrees, $x(t)$ is the position along the $x$-axis and $y(t)$ is the position along the $y$-axis.
 
The third part, $F_\varepsilon\left(\bm{x}\left(t\right),\bm{\omega}\right)$ accounts for the height $h$ in meters above the ocean floor of the sensor attached to the autonomous vehicle, and is  given by
\begin{equation}
F_\varepsilon\left(\bm{x}\left(t\right),\bm{\omega}\right):=\frac{1}{1+e^{p_\epsilon\left(\varepsilon_\text{DE}-\frac{\varepsilon_\text{FOV}}{2} - \varepsilon^b\left(\bm{x}\left(t\right),\bm{\omega}\right)\right)}} + \frac{1}{e^{p_\epsilon\left(\varepsilon^b\left(\bm{x}\left(t\right),\bm{\omega}\right)-\varepsilon_\text{DE}-\frac{\varepsilon_\text{FOV}}{2}\right)}} - 1,
\end{equation}
where $p_\varepsilon$ is a parameter used to adjust the slope of the sigmoidal curves, and
\begin{equation}
\varepsilon^b\left(\bm{x}\left(t\right),\bm{\omega}\right) := \arctan\left(\frac{-h}{\norm{\bm{\omega} - \bm{x}\left(t\right)} }\right),
\end{equation}
where $\varepsilon_\text{FOV}$ and $\varepsilon_\text{DE}$ respectively stand for the vertical FOV, and the downward elevation angle such that the sensor can ensonify the sea floor. For a more thorough description we refer to \cite{9140316}.
\subsection{Trajectories}

In \cite{9140316}, the optimization problem was formulated such that the residual MCM risk is minimized for a given mission time. The residual MCM risk is defined as the probability that a team of autonomous vehicles fails to detect the mines in a search area by the end of an MCM operation. The residual MCM risk is sometimes also referred to as the probability of non-detection.

 In this work,  however we formulated the optimization problem such that the  mission time $T_f$ needed to survey a designated zone $\Omega$ is minimized for a given residual MCM risk,
\begin{equation}
\text{min}\, T_f,
\end{equation}
subjected to
\begin{equation}
 \E[q\left(T_F\right)] :=  \int_\Omega \text{e}^{-\int_0^{T_F} \gamma\left(\bm{x}\left(\tau\right),\bm{\omega}\right)\, d\,\tau}\phi\left(\bm{\omega}\right) d\,\bm{\omega} \leq  \text{residual MCM risk}
\label{eq:exp}
\end{equation}
where $\phi(\cdot)$ is the CDF of the  distribution against which we integrate.  In this case we consider $\phi(\cdot)$ to be the CDF of the Uniform distribution. We note that in our current approach, a Monte Carlo integration scheme is used to compute the expected value of Eq.\,\eqref{eq:exp}.

The position of the autonomous vehicle is given by 

\begin{equation}
\bm{x}\left(t\right) := f(x(t), y(t), \psi(t), r(t)),
\label{eq:pos}
\end{equation}
with $r(t)$ being the turn rate in degrees per second.
The position at time $t$, see Eq.\,\eqref{eq:pos}, is governed by the following differential equations,
\begin{equation}
\begin{aligned}
&\frac{d\,x(t)}{d\,t} = V \text{cos}(\psi(t)) \\
&\frac{d\,y(t)}{d\,t} = V \text{sin}(\psi(t)) \\
&\frac{d\,\psi(t)}{d\,t} = r(t) \\
&\frac{d\,r(t)}{d\,t} = \frac{1}{T}\left(K p(t) - r(t)\right),
\end{aligned}
\end{equation}
where $V$ stands for the speed in $m/s$, $K$ is the Nomoto gain constant with units $s^{-1}$, $T$ is the Nomoto time constant with units ${s}$, and $p(t)$ is the rudder deflection angle given in degrees. The optimization software will control the values for the rudder deflection angle in order to compute a trajectory satisfying the target function and the constraints.

\subsection{Modeling of Ripples}
The presence of sand ripples in the domain is taken into account by multiplying the gamma function, $\gamma\left(\bm{x}\left(\tau\right),\bm{\omega}\right)$ of Eq.\,\eqref{eq:exp},  with a set of functions dividing the rectangular domain $\Omega$ into a zone containing ripples and one not containing ripples,
\begin{equation}
\gamma\left(\bm{x}\left(t\right),\bm{\omega}\right)_{\text{new}} = \gamma\left(\bm{x}\left(t\right),\bm{\omega}\right) \, \text{Dom}\left(\omega_x,\omega_y,\psi(t)\right).
\end{equation}

%

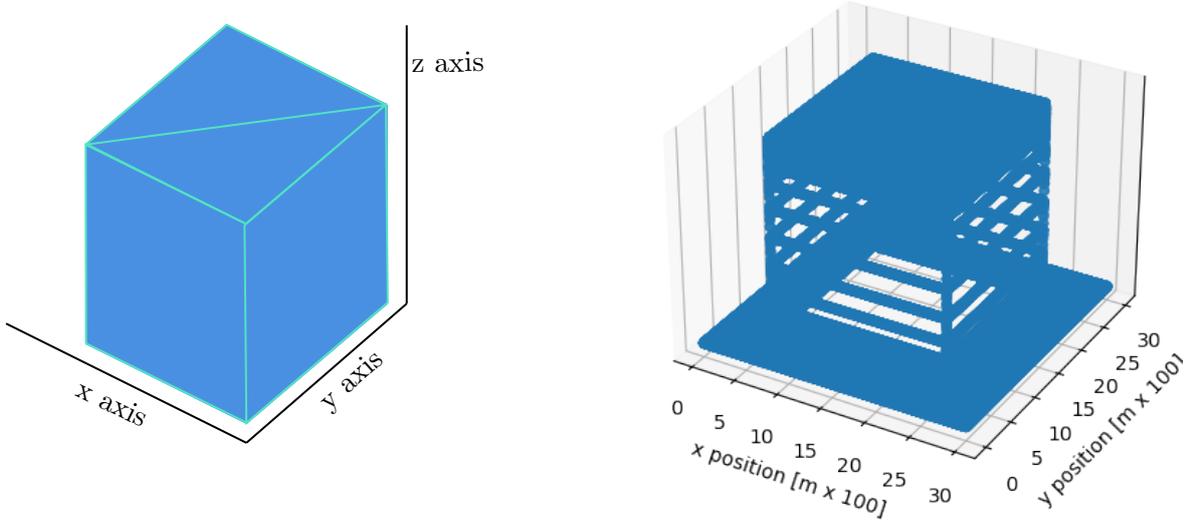
\begin{figure}
  \begin{minipage}{0.5\textwidth}
\tikzset{every picture/.style={line width=0.75pt}} 

\begin{tikzpicture}[x=0.75pt,y=0.75pt,yscale=-1,xscale=1]

\draw    (80,170) -- (200,230) ;
\draw    (200,230) -- (280,160) ;
\draw [color={rgb, 255:red, 74; green, 144; blue, 226 }  ,draw opacity=1 ]   (120,180) -- (200,220) ;
\draw [color={rgb, 255:red, 74; green, 144; blue, 226 }  ,draw opacity=1 ]   (200,220) -- (270,160) ;
\draw    (280,20) -- (280,160) ;
\draw [color={rgb, 255:red, 80; green, 227; blue, 194 }  ,draw opacity=1 ]   (120,80) -- (120,180) ;
\draw [color={rgb, 255:red, 74; green, 144; blue, 226 }  ,draw opacity=1 ]   (270,60) -- (270,160) ;
\draw [color={rgb, 255:red, 74; green, 144; blue, 226 }  ,draw opacity=1 ]   (120,80) -- (200,120) ;
\draw [color={rgb, 255:red, 74; green, 144; blue, 226 }  ,draw opacity=1 ][fill={rgb, 255:red, 74; green, 144; blue, 226 }  ,fill opacity=1 ]   (200,120) -- (200,220) ;
\draw [color={rgb, 255:red, 74; green, 144; blue, 226 }  ,draw opacity=1 ]   (200,120) -- (270,60) ;
\draw [color={rgb, 255:red, 74; green, 144; blue, 226 }  ,draw opacity=1 ]   (190,20) -- (120,80) ;
\draw [color={rgb, 255:red, 74; green, 144; blue, 226 }  ,draw opacity=1 ][fill={rgb, 255:red, 74; green, 144; blue, 226 }  ,fill opacity=1 ]   (190,20) -- (270,60) ;
\draw  [color={rgb, 255:red, 80; green, 227; blue, 194 }  ,draw opacity=1 ][fill={rgb, 255:red, 74; green, 144; blue, 226 }  ,fill opacity=1 ] (190.21,19.77) -- (269.57,59.99) -- (199.55,119.93) -- (120.19,79.7) -- cycle ;
\draw  [color={rgb, 255:red, 80; green, 227; blue, 194 }  ,draw opacity=1 ][fill={rgb, 255:red, 74; green, 144; blue, 226 }  ,fill opacity=1 ] (119.8,79.99) -- (199.88,119.85) -- (200.08,219.86) -- (120,180) -- cycle ;
\draw [color={rgb, 255:red, 80; green, 227; blue, 194 }  ,draw opacity=1 ]   (119.8,79.99) -- (200,120) ;
\draw  [color={rgb, 255:red, 80; green, 227; blue, 194 }  ,draw opacity=1 ][fill={rgb, 255:red, 74; green, 144; blue, 226 }  ,fill opacity=1 ] (199.06,120.57) -- (269.57,59.99) -- (270.5,159.42) -- (199.99,220) -- cycle ;
\draw [color={rgb, 255:red, 80; green, 227; blue, 194 }  ,draw opacity=1 ]   (119.8,79.99) -- (269.57,59.99) ;

\draw (117.8,195.05) node [anchor=north west][inner sep=0.75pt]  [rotate=-26.16] [align=left] {x axis};
\draw (232.42,206.49) node [anchor=north west][inner sep=0.75pt]  [rotate=-318.22] [align=left] {y axis};
\draw (280.82,32.43) node [anchor=north west][inner sep=0.75pt]  [rotate=-358.76] [align=left] {z axis};

\end{tikzpicture}
 \end{minipage} 
 \begin{minipage}{0.5\textwidth}             
\includegraphics[scale=0.75]{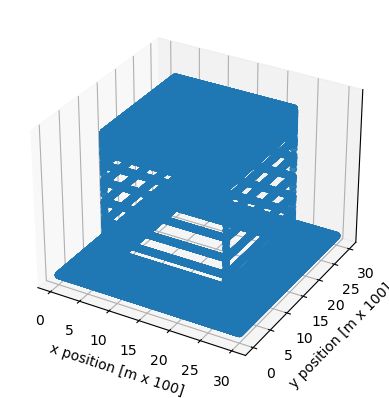}
 \end{minipage}  
 \caption{Schematic view the two-dimensional rectangular function (left) and two dimensional rectangular function on the domain $\Omega = \left[5,25\right]^2$ (right).}
\label{fig:rect_fun}
\end{figure}

\begin{figure}
\centering
\includegraphics[scale=0.5]{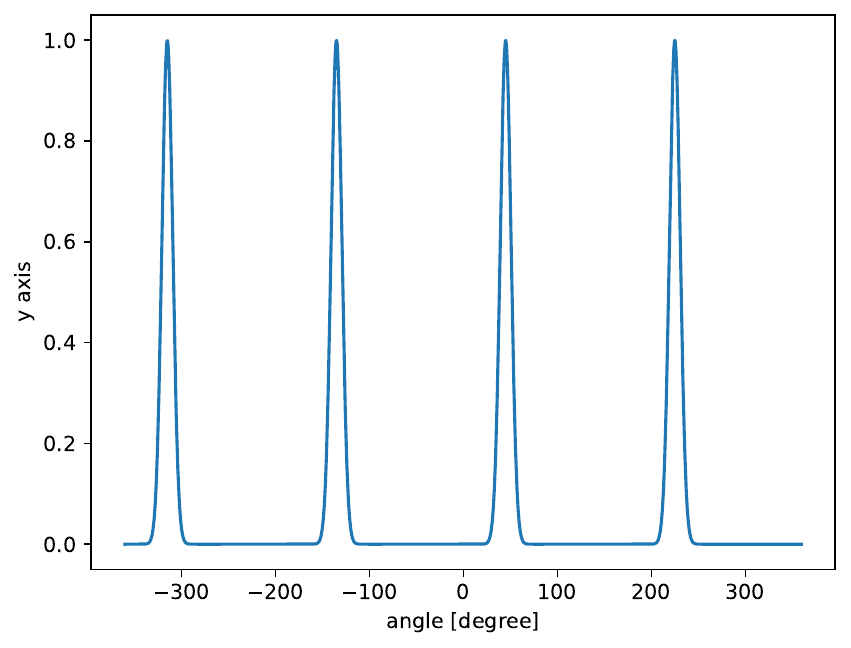}
\caption{The ripple function for ripples at $135^\text{o}$.}
\label{fig:angle}
\end{figure}

The domain function, $\text{Dom}\left(x,y,\theta\right)$, is a combination of the soft rectangular function taken from \cite{SCHAEFER2017368},  see Fig.\,\ref{fig:rect_fun},  and the ripple function see Fig.\,\ref{fig:angle}. The ripple function presented here, models the sand ripples under an angle of $135^\text{o}$. We have chosen that the soft rectangular function divides the square domain $\Omega$ into two triangular zones; one zone where we model sand ripples and one where we do not model sand ripples.
This rectangular function is given as
\begin{equation}
\begin{split}
\text{Rec}\left(x,y\right) := \\ &  \left[\left(\frac{\tanh\left(30\left(x-5\right)\right)-\tanh\left(30\left(x-25\right)\right)}{2}\right)  \right.
  \left.\left(\frac{\tanh\left(30\left(y-x\right)\right)-\tanh\left(30\left(y-25\right)\right)}{2} \right)\right]+ \\
 & \left[\left(\frac{\tanh\left(30\left(x-5\right)\right)-\tanh\left(30\left(x-25\right)\right)}{2}\right) \right. 
  \left.\left(\frac{\tanh\left(30\left(y-25\right)\right)-\tanh\left(30\left(y-x\right)\right)}{2}\right) \right].
 \end{split}
 \label{eq:split}
\end{equation}
The ripple function is given as
\begin{equation}
\text{Ripple}\left(\theta\right) := \sum_{k=-1}^2\exp\left(-\frac{\left(\frac{\theta - \pi / 4 + \pi k}{0.1}\right)^2}{2}\right).
\label{eq:rip}
\end{equation}
 The $\text{Dom}\left(x,y,\theta\right)$ function is given as 
\begin{equation}
\begin{split}
  \text{Dom}\left(x,y,\theta\right) := \\ &   \left[\text{Ripple}\left(\theta\right) \left(\frac{\tanh\left(30\left(x-5\right)\right)-\tanh\left(30\left(x-25\right)\right)}{2}\right.\right) \\& \left(\left.\frac{\tanh\left(30\left(y-x\right)\right)-\tanh\left(30\left(y-25\right)\right)}{2} \right)\right] +  \\ & 
 \left[\left(\frac{\tanh\left(30\left(x-5\right)\right)-\tanh\left(30\left(x-25\right)\right)}{2}\right)  \right. \\&  \left.\left(\frac{\tanh\left(30\left(y-25\right)\right)-\tanh\left(30\left(y-x\right)\right)}{2} \right) \right].
 \end{split}
\end{equation}
The approach as described here will force the optimization software to select trajectories that are perpendicular  to the position of the sand ripples. This is because a contribution to the expected value defined in Eq.\,\eqref{eq:exp} will only be made when the autonomous vehicle is travelling in a perpendicular way to the sand ripples.
%

\section{RESULTS}

\begin{figure}
\centering
\scalebox{0.9}{
\tikzset{every picture/.style={line width=0.75pt}} 

\begin{tikzpicture}[x=0.75pt,y=0.75pt,yscale=-1,xscale=1]
\draw   (110,50) -- (310,50) -- (310,250) -- (110,250) -- cycle ;
\draw    (110,250) -- (310,50) ;
\draw [color={rgb, 255:red, 208; green, 2; blue, 27 }  ,draw opacity=1 ] [dash pattern={on 0.84pt off 2.51pt}]  (110,50) -- (210,150) ;
\draw [color={rgb, 255:red, 208; green, 2; blue, 27 }  ,draw opacity=1 ] [dash pattern={on 0.84pt off 2.51pt}]  (110,70) -- (200,160) ;
\draw [color={rgb, 255:red, 208; green, 2; blue, 27 }  ,draw opacity=1 ] [dash pattern={on 0.84pt off 2.51pt}]  (130,50) -- (220,140) ;
\draw [color={rgb, 255:red, 208; green, 2; blue, 27 }  ,draw opacity=1 ] [dash pattern={on 0.84pt off 2.51pt}]  (150,50) -- (230,130) ;
\draw [color={rgb, 255:red, 208; green, 2; blue, 27 }  ,draw opacity=1 ] [dash pattern={on 0.84pt off 2.51pt}]  (170,50) -- (240,120) ;
\draw [color={rgb, 255:red, 208; green, 2; blue, 27 }  ,draw opacity=1 ] [dash pattern={on 0.84pt off 2.51pt}]  (190,50) -- (250,110) ;
\draw [color={rgb, 255:red, 208; green, 2; blue, 27 }  ,draw opacity=1 ] [dash pattern={on 0.84pt off 2.51pt}]  (210,50) -- (260,100) ;
\draw [color={rgb, 255:red, 208; green, 2; blue, 27 }  ,draw opacity=1 ] [dash pattern={on 0.84pt off 2.51pt}]  (230,50) -- (270,90) ;
\draw [color={rgb, 255:red, 208; green, 2; blue, 27 }  ,draw opacity=1 ] [dash pattern={on 0.84pt off 2.51pt}]  (250,50) -- (280,80) ;
\draw [color={rgb, 255:red, 208; green, 2; blue, 27 }  ,draw opacity=1 ] [dash pattern={on 0.84pt off 2.51pt}]  (270,50) -- (290,70) ;
\draw [color={rgb, 255:red, 208; green, 2; blue, 27 }  ,draw opacity=1 ] [dash pattern={on 0.84pt off 2.51pt}]  (290,50) -- (300,60) ;
\draw [color={rgb, 255:red, 208; green, 2; blue, 27 }  ,draw opacity=1 ] [dash pattern={on 0.84pt off 2.51pt}]  (110,90) -- (190,170) ;
\draw [color={rgb, 255:red, 208; green, 2; blue, 27 }  ,draw opacity=1 ] [dash pattern={on 0.84pt off 2.51pt}]  (110,110) -- (180,180) ;
\draw [color={rgb, 255:red, 208; green, 2; blue, 27 }  ,draw opacity=1 ] [dash pattern={on 0.84pt off 2.51pt}]  (110,130) -- (170,190) ;
\draw [color={rgb, 255:red, 208; green, 2; blue, 27 }  ,draw opacity=1 ] [dash pattern={on 0.84pt off 2.51pt}]  (110,150) -- (160,200) ;
\draw [color={rgb, 255:red, 208; green, 2; blue, 27 }  ,draw opacity=1 ] [dash pattern={on 0.84pt off 2.51pt}]  (110,170) -- (150,210) ;
\draw [color={rgb, 255:red, 208; green, 2; blue, 27 }  ,draw opacity=1 ] [dash pattern={on 0.84pt off 2.51pt}]  (110,190) -- (140,220) ;
\draw [color={rgb, 255:red, 208; green, 2; blue, 27 }  ,draw opacity=1 ] [dash pattern={on 0.84pt off 2.51pt}]  (110,210) -- (130,230) ;
\draw [color={rgb, 255:red, 208; green, 2; blue, 27 }  ,draw opacity=1 ] [dash pattern={on 0.84pt off 2.51pt}]  (110,230) -- (120,240) ;
\end{tikzpicture}}
\caption{Delimitation of the ripple zone. Sand ripples are present in the upper-left triangle at $135^\text{o}$ and depicted by the red dotted line.}
\label{fig:zone}
\end{figure}
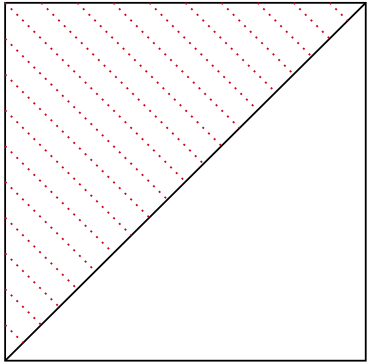

In this section we present the results when considering one and two autonomous vehicles in a square domain $\Omega = \left[5,25\right]^2$ where sand ripples are present in the upper-left triangle at $135^\text{o}$, see Fig.\,\ref{fig:zone}. We consider three different cases. The first case consists of computing the trajectory when no ripples are present in the zone. In the second case, we do not recompute the trajectory, i.e., we use the same trajectory as in case 1, but take into account the presence of sand ripples. This second case mainly serves to illustrate the impact on the residual MCM risk when not taking into account the presence of sand ripples when these are present. The third case shows an updated trajectory which accounts for the presence of the ripples. A brief description of the cases is also given in Tab.\,\ref{Tab:cases}. For Case 1 and case 3, we requested a residual MCM risk of $10\,\%$. All simulations have been performed on a 16  core Intel i7-12850HX processor with 32GB RAM. The different parameters for our simulations are shown in Tab.\,\ref{Tab:Param}.

\begin{table}[h]
\centering
\scalebox{1.0}{
\begin{tabular}{cc}
\toprule
 {Parameter name} &  {Value} \\
 \cmidrule(rl{4pt}){1-2}  
  $\alpha_{\text{FOV}}$ &  $120.0^\text{o}$  \\
 $h$ &  20.0\,$m$ \\
 $\sigma$ &  9.0\,[/]  \\
 $\lambda$ & 20.0\,$s^{-1}$\\
 $a$ & 5.2\,$dB/km$ \\
 $\varepsilon_{\text{FOV}}$ & $5.0^\text{o}$\\
  $\varepsilon_{\text{DE}}$ & $-6.0^\text{o}$\\
 $\text{FOM}$ & 72.0\,[/]\\
 $p_\alpha$ & 25.0\,[/]\\
 $p_\varepsilon$ & 400.0\,[/]\\
 $\text{V}$ & 2.5\,$m/s$\\
  $\text{T}$ & 0.5\,$s$\\
 $\text{K}$ & 5.0\,$s^{-1}$\\
\bottomrule
\end{tabular}}
\caption{Parameters for the simulations.}
\label{Tab:Param}
\end{table}

\begin{table}
\centering
\scalebox{1.0}{
\begin{tabular}{cc}
\toprule
 {Case number} &  {Description} \\
 \cmidrule(rl{4pt}){1-2}  
  1 &  Computation of trajectory with no ripples present  \\
 2 &  Identical trajectory as case 1 but with ripples  present  \\
 3 &  Computation of an updated trajectory accounting for the presence of ripples    \\
\bottomrule
\end{tabular}}
\caption{Description of the different cases.}
\label{Tab:cases}
\end{table} 

\subsection{One Autonomous Vehicle}
\begin{figure}
\hspace{1.5cm}
\begin{minipage}{0.45\textwidth}
\scalebox{0.9}{
\begin{tikzpicture}

    \node at (3.5,6.5) {\large \bf---------Case 1---------};

    \node at (3.5,5.5) {\large Residual MCM risk: 10.00\%};

\end{tikzpicture}}
\end{minipage}
\hspace{1.cm}
\begin{minipage}{0.45\textwidth}
\scalebox{0.9}{
\begin{tikzpicture}

    \node at (3.5,6.5) {\large \bf---------Case 2---------};

    \node at (3.5,5.5) {\large Residual MCM risk: 55.75\%};

\end{tikzpicture}}
\end{minipage}
\\
  \begin{minipage}{0.55\textwidth}   
\includegraphics[scale=0.4]{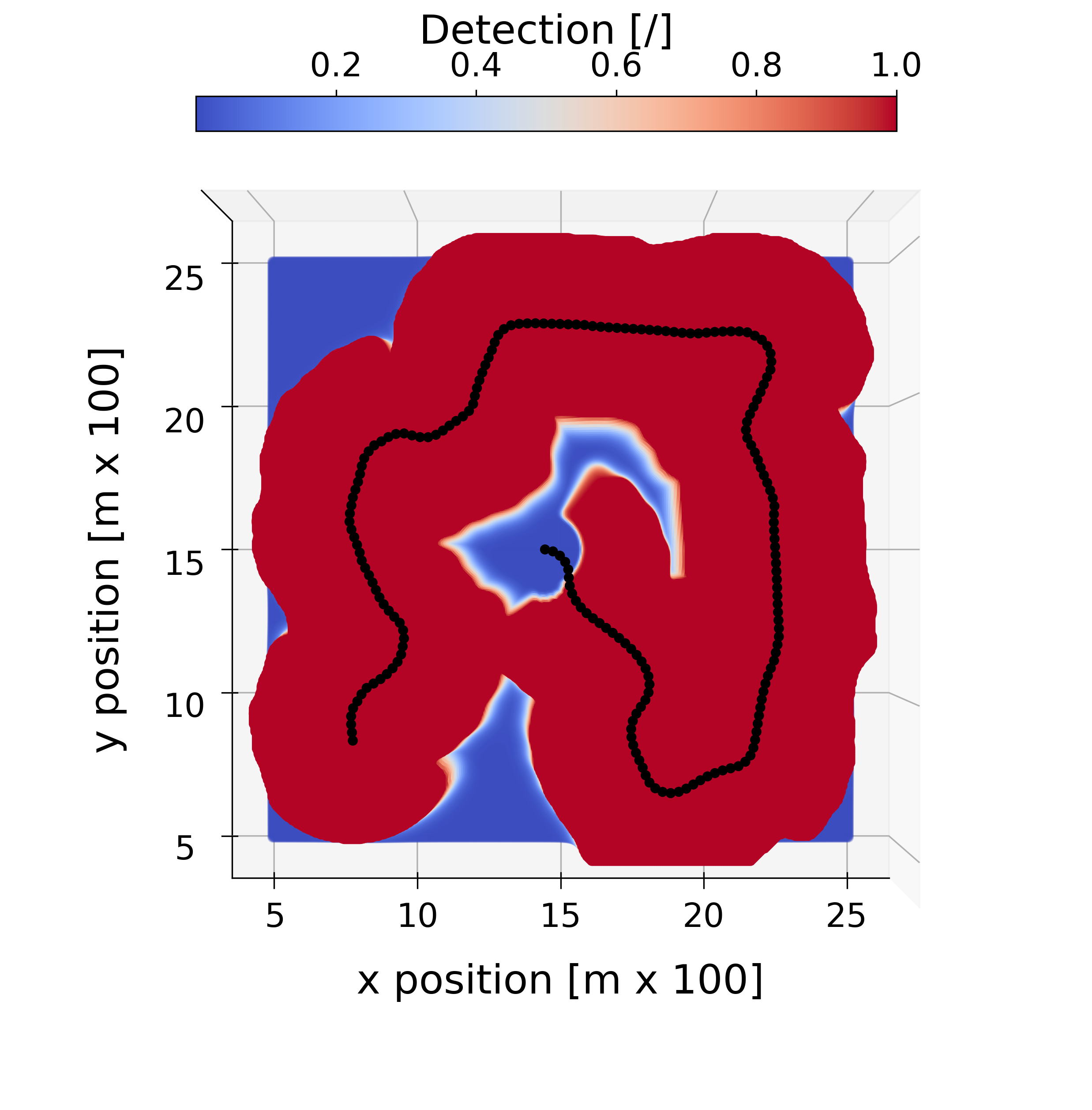}
 \end{minipage}  
  \hspace{-0.4cm}
  \begin{minipage}{0.5\textwidth}
\includegraphics[scale=0.4]{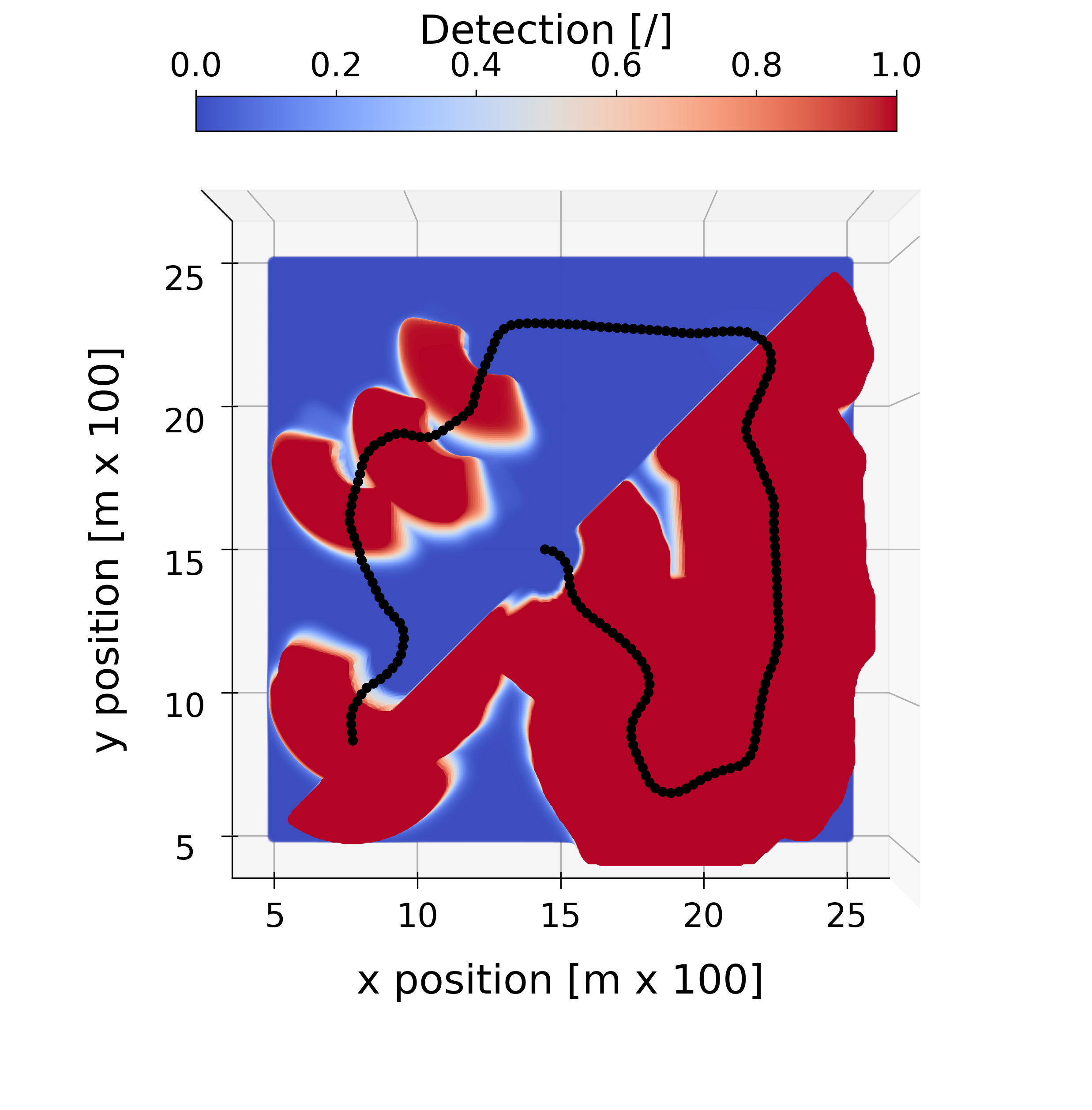}
 \end{minipage} 
 \\
\begin{minipage}{0.45\textwidth}
\vspace{1.cm}
\hspace{5.5cm}
\scalebox{0.9}{
\begin{tikzpicture}

    \node at (3.5,6.5) {\large \bf---------Case 3---------};
    
        \node at (3.5,5.5) {\large Residual MCM risk: 10.00\%};

\end{tikzpicture}}
\end{minipage}
\vspace{0.25cm}
\\
  \begin{minipage}{1.0\textwidth}
    \hspace{4.cm}
\includegraphics[scale=0.4]{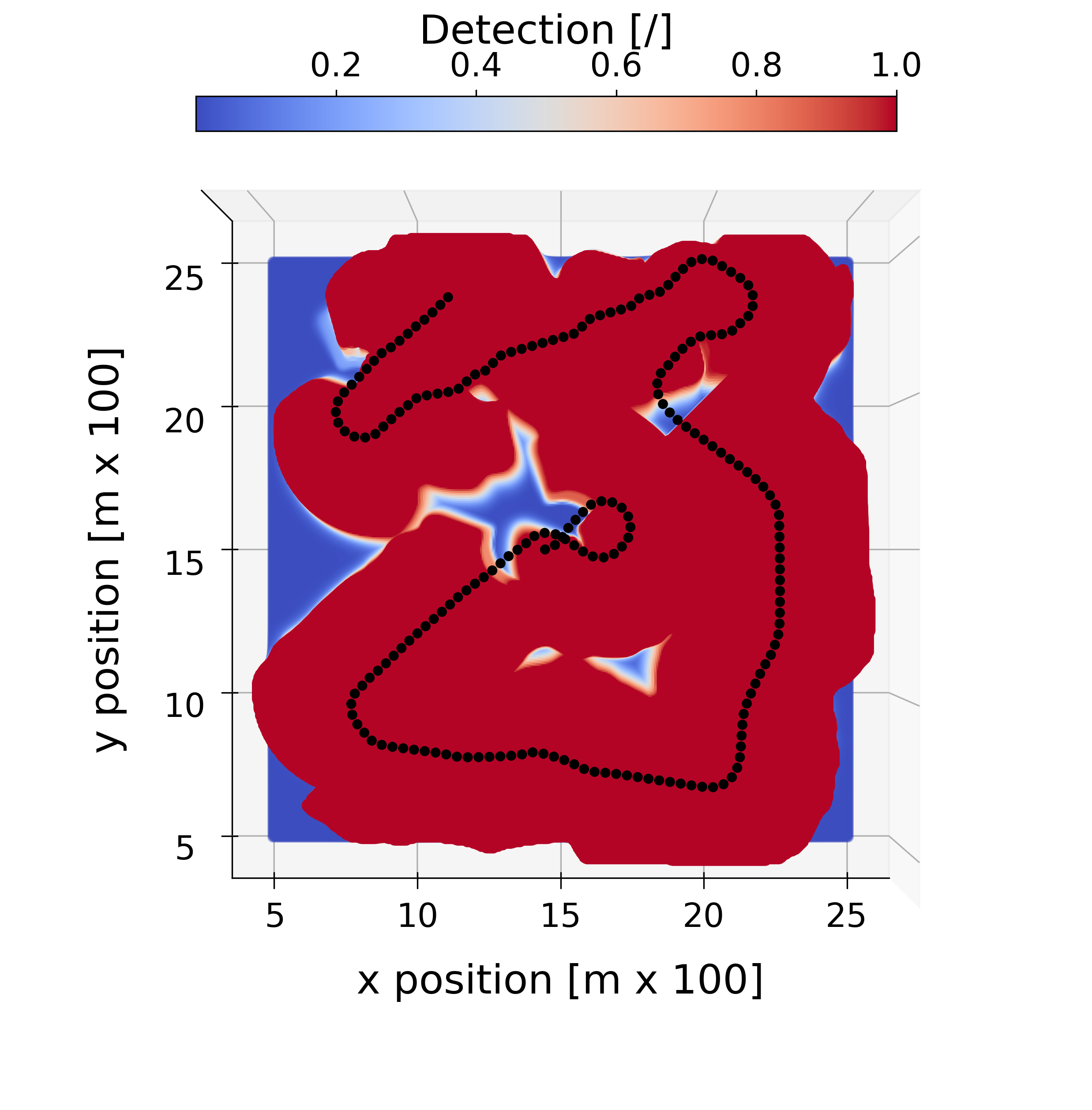}
 \end{minipage} 
 \caption{Trajectories and residual MCM risks, see Tab.\,\ref{Tab:res_1S}, for the results of the  different cases as defined in Tab.\,\ref{Tab:cases}. The red color represents the zone that has been surveyed by the sensor, while the blue zone has not been surveyed. The detection probability of a mine by the sensor in the red zone is 1.0 or 100\,\% and 0\,\% in the blue zone. The black line represents the trajectory of the autonomous vehicle.}
\label{fig:res_1s}
\end{figure}
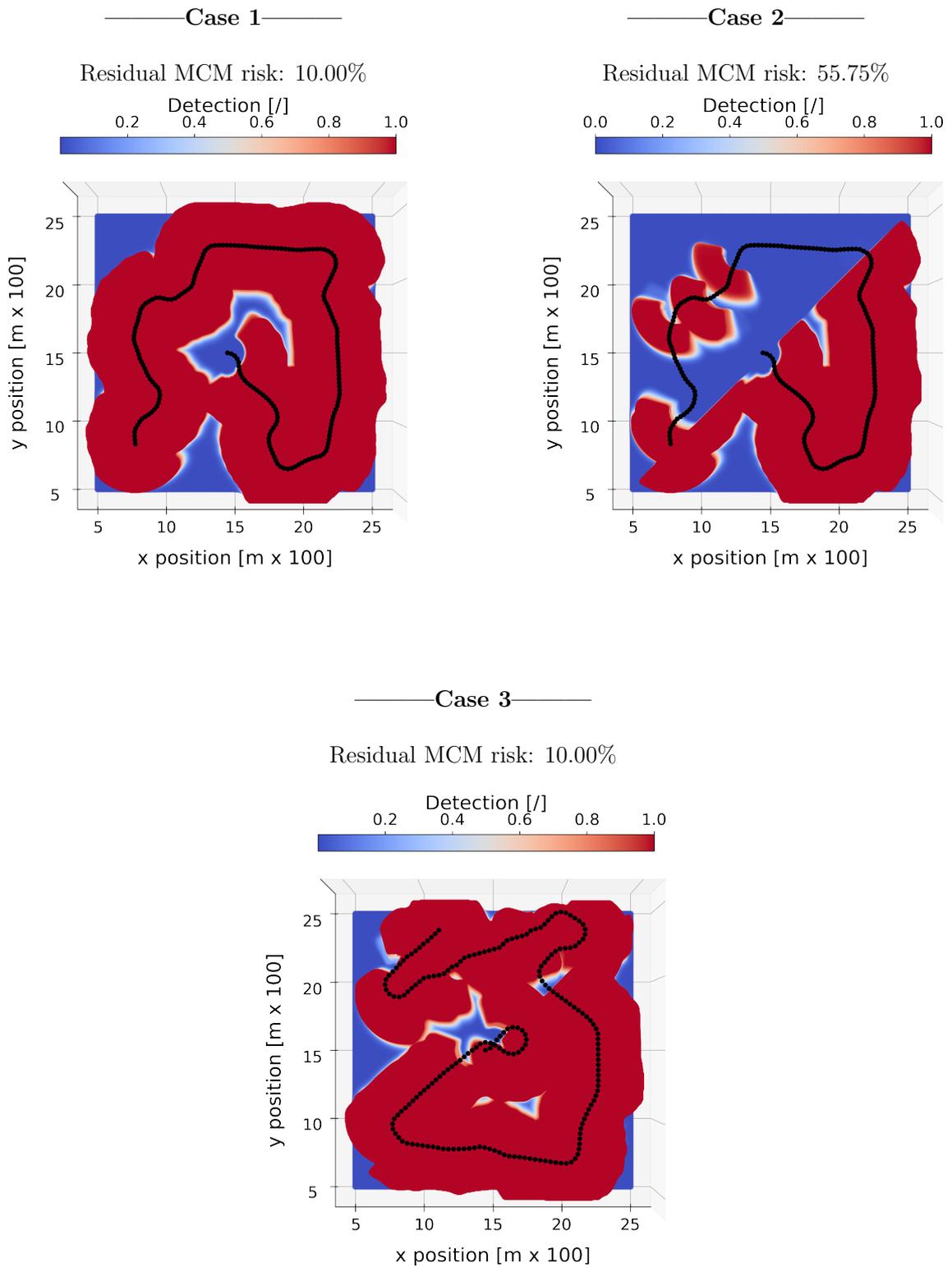

We present the results for the different cases when considering only one autonomous vehicle, in Tab.\,\ref{Tab:res_1S} and Fig.\,\ref{fig:res_1s}. In Fig.\,\ref{fig:res_1s}, the red color represents the zone that has been surveyed by the sensor, while the blue zone has not been surveyed. The detection probability of a mine by the sensor in the red zone is 1.0 or 100\,\% and 0\,\% in the blue zone. The black line represents the trajectory of the autonomous vehicle. The coordinates of the  starting point for the autonomous vehicle for all the cases is $\left(14.5,15.0\right)$. For case 1, the residual MCM risk is  $10.00$\,\%. The mission time necessary to achieve this percentage is 2088.00 seconds. For case 2, the residual MCM risk is  $55.75$\,\%, i.e., an increase in risk of  $45.75$\,\% point. This result shows the impact on the residual MCM risk when not taking the ripples into account when these are present.  For case 3, the residual MCM risk is again $10.00$\,\% but with a longer mission time, i.e., 2778.85 seconds instead of 2088.00 seconds when compared to case 1. Accounting for sand ripples thus increases the  mission duration time. The total computation time needed to calculate the trajectory also increases when accounting for sand ripples, from 20.67 seconds to 93.18 seconds.

\begin{table}
\centering
\scalebox{1.0}{
\begin{tabular}{ccccc}
\toprule
 {Case} &  {Residual MCM risk} &   {Mission time} & {Computation time} \\
 \cmidrule(rl{4pt}){1-4}  
 1 &  10.00\,\% &  2088.00\,seconds &20.67\,seconds      \\
 2 &  55.75\,\% &  2088.00\,seconds &/     \\
 3 &  10.00\,\% &  2778.85\,seconds &93.18\,seconds       \\
\bottomrule
\end{tabular}}
\caption{Results when considering one autonomous vehicle.}
\label{Tab:res_1S}
\end{table} 

\subsection{Two Autonomous Vehicles}

\begin{figure}
\hspace{1.5cm}
\begin{minipage}{0.45\textwidth}
\scalebox{0.9}{
\begin{tikzpicture}

    \node at (3.5,6.5) {\large \bf---------Case 1---------};

    \node at (3.5,5.5) {\large Residual MCM risk: 10.00\%};

\end{tikzpicture}}
\end{minipage}
\hspace{1.cm}
\begin{minipage}{0.45\textwidth}
\scalebox{0.9}{
\begin{tikzpicture}

    \node at (3.5,6.5) {\large \bf---------Case 2---------};

    \node at (3.5,5.5) {\large Residual MCM risk: 57.18\%};

\end{tikzpicture}}
\end{minipage}
\\
  \begin{minipage}{0.55\textwidth}             
\includegraphics[scale=0.4]{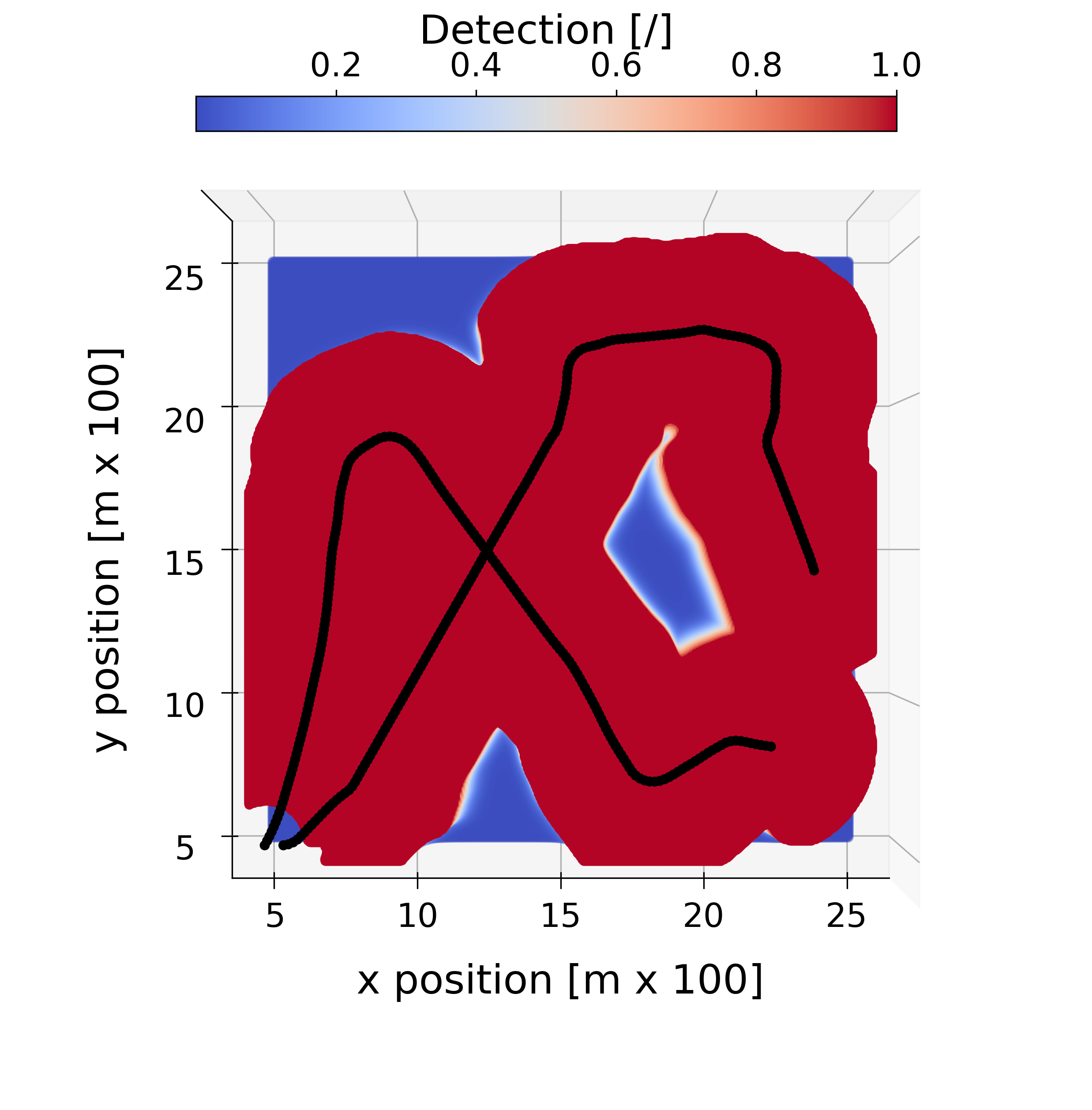}
 \end{minipage}  
  \hspace{-0.4cm}
  \begin{minipage}{0.5\textwidth}
\includegraphics[scale=0.4]{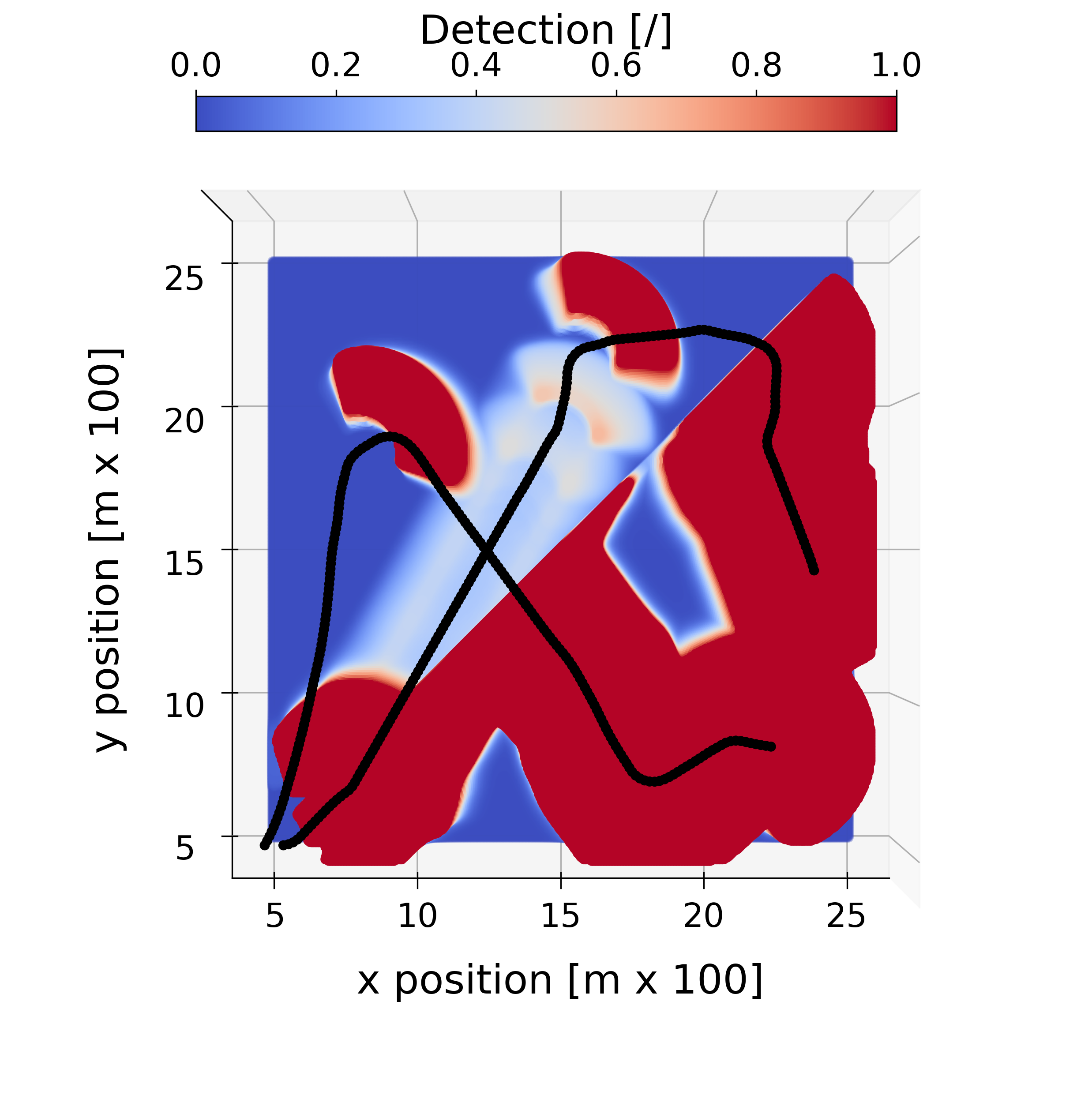}
 \end{minipage} 
 \\
\begin{minipage}{0.45\textwidth}
\vspace{1.cm}
\hspace{5.5cm}
\scalebox{0.9}{
\begin{tikzpicture}

    \node at (3.5,6.5) {\large \bf---------Case 3---------};
    
        \node at (3.5,5.5) {\large Residual MCM risk: 10.00\%};

\end{tikzpicture}}
\end{minipage}
\vspace{0.25cm}
\\
  \begin{minipage}{1.0\textwidth}
    \hspace{4.cm}
\includegraphics[scale=0.4]{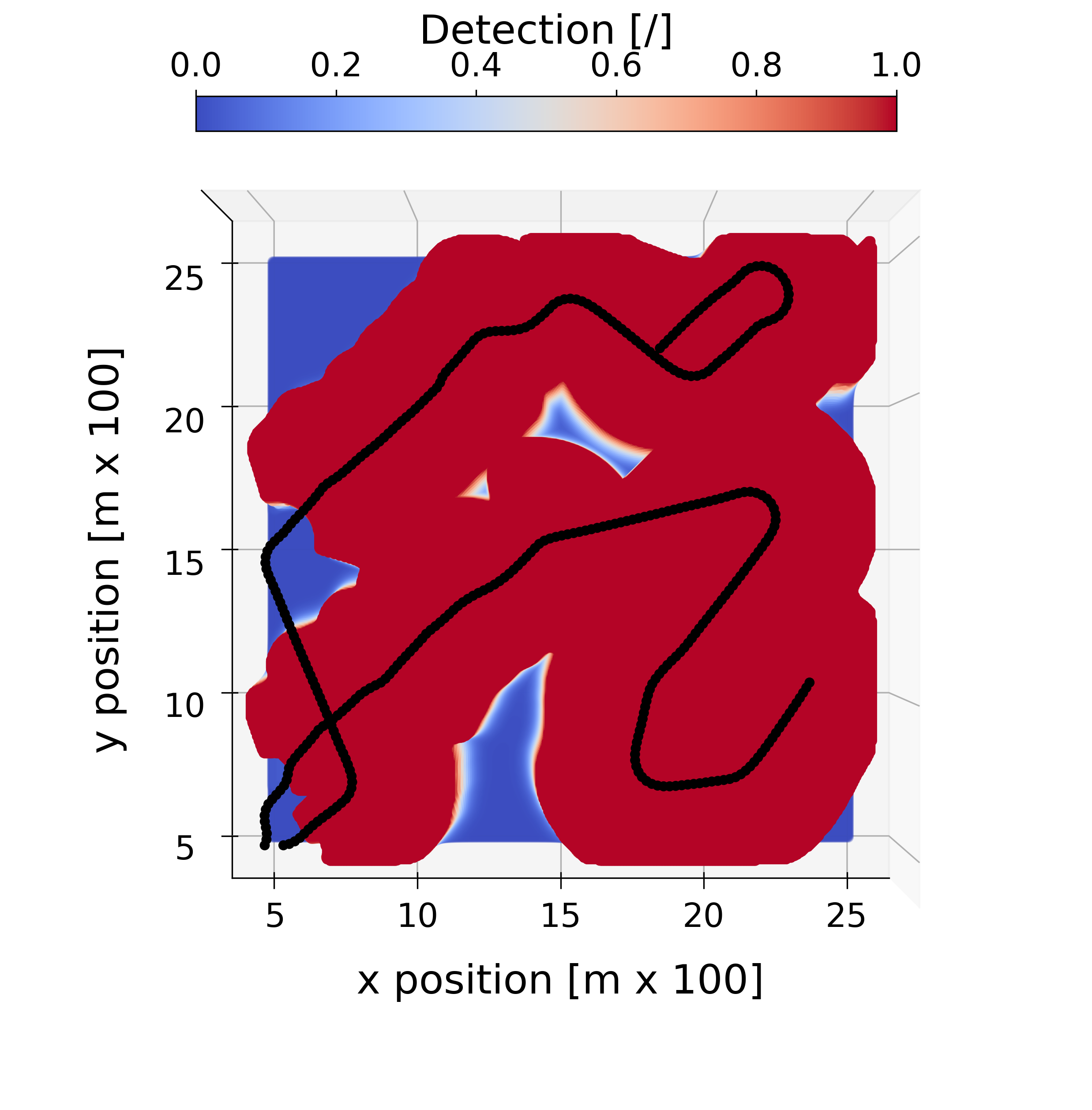}
 \end{minipage} 
 \caption{Trajectories and residual MCM risks, see Tab.\,\ref{Tab:res_2S}, for the results of the different cases as defined in Tab.\,\ref{Tab:cases}. The red color represents the zone that has been surveyed by the combined sensors, while the blue zone has not been surveyed. The detection probability of a mine by the combined sensors in the red zone is 1.0 or 100\,\% and 0\,\% in the blue zone. The black lines represent the trajectories of the autonomous vehicles.}

 \label{fig:res_2s}
\end{figure}
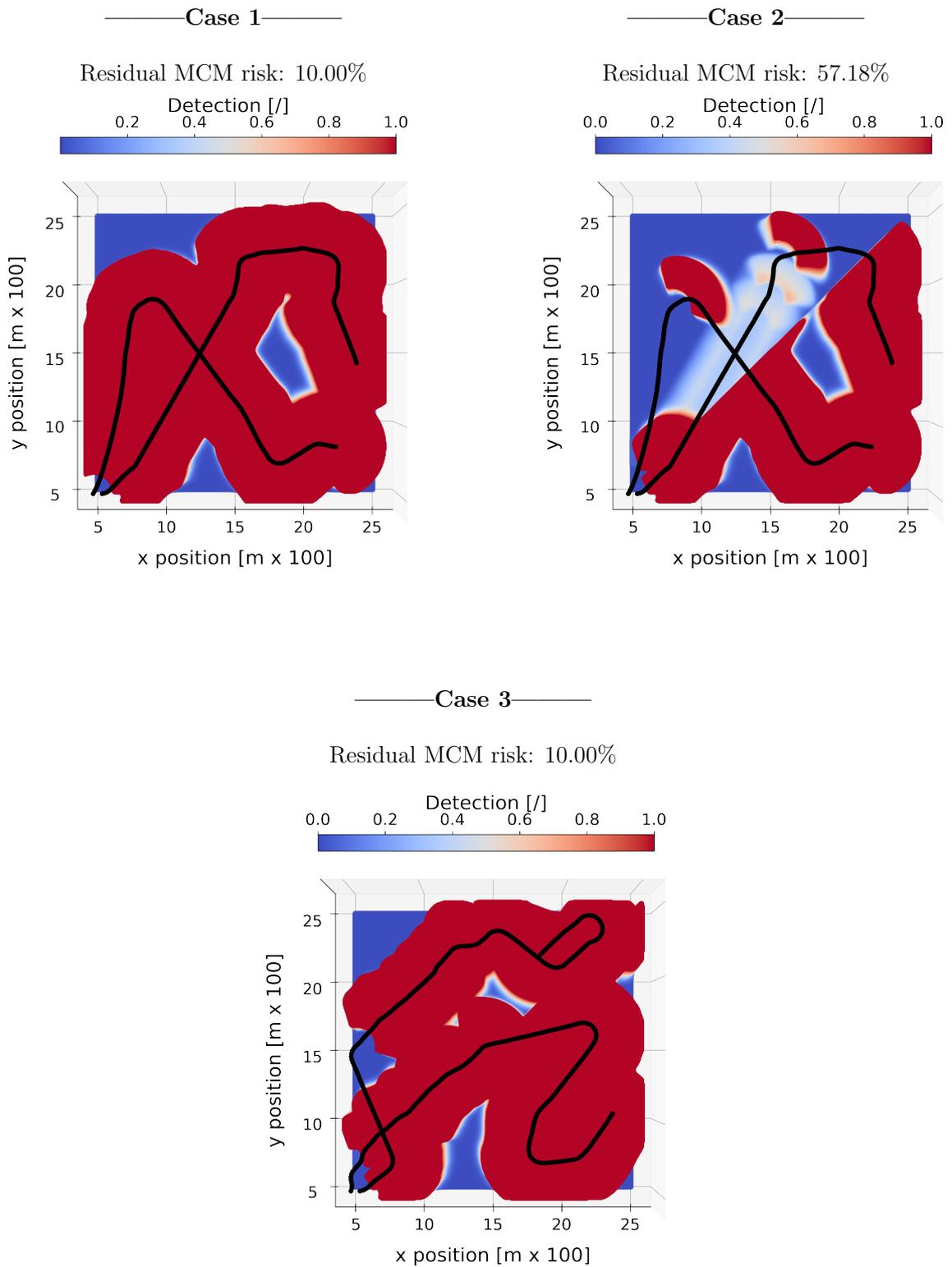

We now present the results for the different cases when considering two autonomous vehicles, in Tab.\,\ref{Tab:res_2S} and Fig.\,\ref{fig:res_2s}. In Fig.\,\ref{fig:res_2s}, the red color represents the zone that has been surveyed by the combined sensors, while the blue zone has not been surveyed. The detection probability of a mine by the combined sensors in the red zone is 1.0 or 100\,\% and 0\,\% in the blue zone. The black lines represent the trajectories of the autonomous vehicles. The conclusions  are similar to the ones made when considering only one autonomous vehicle. For case 1, the residual MCM risk is 10.00\,\%. The mission time necessary to achieve this percentage is 1303.88\,seconds. For case 2, the residual MCM risk is   57.18\,\%, i.e., an increase in  risk of  47.18\,\% point. For case 3, the residual MCM risk is again 10.00\,\% but with a longer mission time, i.e., 1548.00\,seconds instead of 1303.88\,seconds when compared to case 1. The  mission duration time again increases when considering sand ripples. The total computation time also increases, from 64.22 seconds for case 1 to 118.25 seconds for case 3.

\begin{table}
\centering
\scalebox{1.0}{
\begin{tabular}{ccccc}
\toprule
 {Case} &  {Residual MCM risk} &   {Mission time} & {Computation time} \\
 \cmidrule(rl{4pt}){1-4}  
 1 &  10.00\,\% &  1303.88\,sec &64.22\,sec      \\
 2 &  57.18\,\% &  1303.88\,sec & /     \\
 3 &  10.00\,\% &  1548.00\,sec &118.25\,sec       \\
\bottomrule
\end{tabular}}
\caption{Results when considering two autonomous vehicles.}
\label{Tab:res_2S}
\end{table}

\section{CONCLUSION}
In this work we presented a novel way to account for sand ripples when computing trajectories for autonomous vehicles by means of a stochastic optimal control approach for mine countermeasure operations in a two-dimensional space. In order to model the presence of sand ripples in a given domain, we constructed a function based on the two-dimensional rectangular function, which we multiplied with the function describing the sensor model. This approach forces the optimization software to compute trajectories for autonomous vehicles that are perpendicular to the sand ripples in the zone where  ripples are present. We implemented this for up to two autonomous  vehicles in a two dimensional  domain $\Omega = \left[5,25\right]^2$, and presented results showing the trajectories for a residual MCM risk of  $10\,\%$ when considering sand ripples in the domain. We found that the mission time and computational time are larger when simulating the presence of sand ripples. Future work will focus on investigating how to compute trajectories in a non-square domain. Additionally, we  plan to investigate if a quasi-Monte Carlo integration scheme for computing the expected value of the residual MCM risk integral can yield a speed-up.

\clearpage
\section{REFERENCES} 
\renewcommand{\refname}{} 
\begingroup
\setstretch{0.8}
\bibliographystyle{ieeetr} 
\bibliography{refs} 
\endgroup

\end{document}